\documentclass[12pt]{article}

\usepackage{amsmath,amssymb,amsthm,algorithm,algpseudocode,nicematrix,hyperref}

\newcommand{\cur}{\mathrm{curl}}

\newcommand{\J}{\mathcal{J}}
\newcommand{\T}{\mathcal{T}}
\newcommand{\C}{\mathcal{C}}
\newcommand{\I}{\mathcal{I}}
\newcommand{\R}{\mathbb{R}}
\newcommand{\E}{\mathcal{E}}
\newcommand{\id}{\mathrm{d}}
\DeclareRobustCommand{\Chi}{{\mathpalette\irchi\relax}}
\newcommand{\irchi}[2]{\raisebox{\depth}{$#1\chi$}}

\begin{document}

\setcounter{page}{1}

\title{Electro-thermal topology optimization of an electric machine by the topological derivative considering drive cycles}
\author{Nepomuk~Krenn$^1$         \and
        Théodore~Cherière $^2$ \and
Sebastian~Schöps$^3$ \and
Peter~Gangl$^1$
}
\date{$^1$RICAM, Austrian Academy of Sciences, \\
Altenberger Stra{\ss}e 69, 4040 Linz, Austria\\
$^2$Université Paris-Saclay, Centralesupélec, CNRS, GeePs, \\
91190 Gif-sur-Yvette, France\\
$^3$Computational Electromagentics Group, TU Darmstadt, \\
Schloßgartenstr. 8, 64289 Darmstadt, Germany}

\maketitle
\begin{abstract}
We consider a 2d permanent magnet synchronous machine operating in a sequence of static operating points coming from a drive cycle. We aim to find a rotor design which maximizes the efficiency defined as the quotient of input and output energy considering Joule losses in the stator and eddy current losses in the permanent magnets. A coupled electromagnetic-thermal analysis of the rotor considers the eddy current losses as heat source and adds a temperature constraint to avoid damage of the permanent magnets. Additionally we impose Von-Mises stress constraints to maintain the mechanical integrity of the design. To solve the resulting free form topology optimization problem we use a level set description of the design and the topological derivative as sensitivity information. We show the effect of these constraints at very high speeds which is a trend in recent machine development.

\end{abstract}
\section{Introduction}
Due to the increasing demand for electric machines with high energy density and high efficiency, design optimization became a popular tool to improve the performance of electric machines, especially of permanent magnet synchronous machines (PMSM), which proved to provide the highest power density. A common approach is to optimize a finite number of parameters that describe the geometry, such as positions, widths, lengths, and angles, which also limits the space of possible designs. Besides evolutionary algorithms \cite{Kim2024}, which get expensive for large and more complex design spaces, one can use gradient-based algorithms \cite{Wiesheu2024}, which show a significant speedup of runtime. The efficiency of gradient-based algorithm enables wider optimization spaces that do not rely on parametrization. For instance, another way to optimize the design is to choose an initial topology and optimize the material interfaces, i.e., the shape of the different parts \cite{Kuci2021},\cite{Gangl_Annette}. 

In this work, we focus on topology optimization, which allows for the largest design variations. The most common class of methods, introduced in mechanical engineering \cite{Bendsøe1989}, are density methods. There, one uses a density variable to describe the material distribution. In order to enable gradient computations, one has to introduce material interpolations yielding unphysical intermediate materials. To erase them from the final designs, one applies penalization techniques, as in the widely used SIMP (Solid Isotropic Material Penalization) method, which was recently successfully applied to the multi-material topology optimization of electric machines \cite{Theodore2022}. For more details on topology optimization in electric machines, we refer to the recent review paper \cite{lucchinitorchio_review}.

We will use the approach introduced in \cite{Amstutz_Levelset}, describing the design by a continuous level set function. This has the advantage of crisp material interfaces at any iteration of the optimization. The evolution of this level set function is driven by the topological derivative, the pointwise sensitivity of the objective with respect to material perturbations.

While plenty of the mentioned work focuses on optimizing machines for a single operating point (OP), we will consider a sequence coming from a standardized drive cycle, as it was done, e.g., in \cite{drivecycle_extern}. This is an important step towards real-world applications, assuring the performance of a traction machine under realistic conditions. Following the trend of high-speed machines \cite{Li2016}, we scale the drive cycle up to 27000 rpm, yielding very high efficiencies. This is also our objective: Efficiency defined as the ratio between mechanical output energy and electrical input energy 
\begin{align}\label{eq:efficiency}
    \E(\Omega)=\frac{\int_0^TP^m(\Omega)\id t}{\int_0^TP^e(\Omega)\id t}.
\end{align}
The latter one can be split into the transformed mechanical power and dissipated loss power $P^e(\Omega)=P^m(\Omega)+P^l(\Omega).$ The dominating part of the losses are Joule losses, in which we assume AC losses \cite{Hajji2024} included $P^J(\Omega)=R_SI(\Omega)^2$. We will also consider the eddy current (EC) losses $P^{\mathrm{EC}}(\Omega)$ in the permanent magnets (PMs), which become more and more important with increasing rotation speed $\omega$.

This dissipating power heats up the PMs, which are both very sensitive to high temperatures and difficult to cool, since they may be buried in the rotating rotor. We propose a method, based on ideas from \cite{Amstutz_Constraint} and \cite{AndradeNovotnyLaurain2024}, to constrain the maximum temperature within the PMs integrated in the design optimization. We also include constraints on the Von-Mises stress to ensure mechanical integrity of the design, as it was done in \cite{Amstutz_vonMises},\cite{Holley2023}. Both temperature and mechanical constraint were considered in \cite{Krenn2025}, optimizing a PMSM for a single operating point, which will generalize to drive cycles here.

The rest of this paper is structured as follows: First, we introduce the machine model and the corresponding physical equations for a single, fixed operating point in Section \ref{sec:physics}, which is generalized to drive cycles in Section \ref{sec:drivecycle}. In Section \ref{sec:topop} we introduce the level set method, based on the topological derivative, to optimize the design of an electric machine. Section \ref{sec:constraints} deals with the incorporation of maximum temperature and stress constraints before showing numerical results in Section \ref{sec:results} and concluding in Section \ref{sec:conclusion}.

\section{Physical model}\label{sec:physics}
We first introduce the physical problem in steady state operation for a fixed OP. In the next section we will adapt this to analyze a full drive cycle.
\begin{figure}
    \centering
    \includegraphics[width=0.7\linewidth]{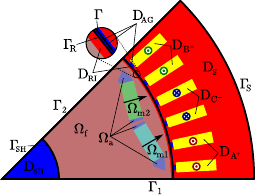}
    \caption{One pole of the machine $D_{\mathrm{all}}$ with outer boundary $\Gamma_S$ and radial boundaries $\Gamma_1, \Gamma_2$, consisting of rotor, stator and airgap. 
    }
    \label{fig:machine}
\end{figure}
\begin{table}
    \centering
    \begin{tabular}{l|l}
         Number of pole pairs $N_{\mathrm{pp}}$ and slots&  4, 48\\
         Inner and outer radius rotor& 17.7mm, 52.4mm  \\
         Inner and outer radius stator& 52.8mm, 77.3mm  \\
         Axial length $\ell_z$&90mm\\
         Stator resistance $R_S$&3.2$\Omega$\\
         Number of turns per slot, fill factor & 60, 0.6\\
         Magnetization direction $\varphi_1,\varphi_2$&  $30^\circ, 15^\circ$\\
         Remanent flux density $B_R$& 1.216T, \\
         Electric conductivity of magnet $\sigma_m$&$6.7\times 10^5$S/m\\
         Thermal conductivities $\lambda_f,\lambda_m,\lambda_a$&16, 9, 0.05 W/(mK)\\
         Thermal robin transfer coefficients $h_{SH},h_{AG}$&0.235, 260 W/(m²K)\\
         Ambient temperature $\theta_0$&$40^\circ$C\\
         Mass densities $\rho_f,\rho_m$&7.65, 8.4\,$\mathrm{g/cm^3}$\\
         Young's modulus $E_f, E_a, E_m$&200, 0.2, 0.2\,GPa\\
         Possion ratio $\nu_f, \nu_a, \nu_m$& 1/3
    \end{tabular}
    \caption{Machine geometry and material data}
    \label{tab:data}
\end{table}
\subsection{Machine model}
We consider a 2d model of one pole of the PMSM $D_{\mathrm{all}}$, depicted in Figure \ref{fig:machine}. 
The rotor splits into three parts: The non-ferromagnetic shaft $D_{\mathrm{SH}}$, a fixed iron ring at the airgap $D_{\mathrm{RI}}$ and the design domain $D$, the latter one consisting of iron in $\Omega_f$, air in $\Omega_a$ and PMs in $\Omega_{m_1}$ and $\Omega_{m_2}$ with fixed remanent flux density $B_R=1.216\mathrm{T}$ but different orientation $\varphi_1,\varphi_2$, respectively. Further machine parameters are given in Table \ref{tab:data}. We denote the actual material configuration in $D$ by
\begin{align*}
     \Omega=(\Omega_f,\Omega_{m_1},\Omega_{m_2},\Omega_a).
\end{align*}
\subsection{Electromagnetic equations}
Electromagnetic phenomena in this machine are described by the equations of 2d nonlinear magnetoquasistatics, a simplification of the Maxwell's equations for low frequency applications. They yield the third component of the magnetic vector potential $a$, a scalar function varying in space and time. After a temporal discretization by $N=11$ equidistant positions on a sixth of an electrical period $t_1=\frac{1}{6N_{\mathrm{pp}}\omega}$, i.e. a mechanical rotation by $15^\circ$, this reads for $n=0,...,N-1$
\begin{align}\label{eq:MQS}
\begin{aligned}
    \frac{\sigma_{\Omega^n}}{\tau}(a^n-a^{n-1})+\cur h_{\Omega^n}(\cur a^n)&=j^n,\\
    a^n|_{\Gamma_S}=0,\quad
    a^n|_{\Gamma_1}&=-a^n|_{\Gamma_2},\\
    \sigma_{\Omega^{-1}}a^{-1}&=\sigma_{\Omega^{N-1}}a^{N-1}.
\end{aligned}
\end{align}
The subindex $._{\Omega^n}$ denotes the dependence of the material laws on the material configuration, which is moving together with the rotor. Each timestep $n$ of size $\tau=t_1N^{-1}$, can be associated to the mechanical rotor angle $15^\circ\frac{n}{N}.$ The last equation leads to a temporally periodic behavior, i.e. steady state operation. We use for iron a nonlinear BH curve taken from \cite{Krenn2025}, 
for non-ferromagnetic material $h_a(b)=\nu_0b$ 
and for PMs $h_{m_i}(b)=\nu_m(b-B_R(\cos\varphi_i,\sin\varphi_i)^T)$, with values given in Table \ref{tab:data}. The machine is excited by a three phase current in distributed winding described by the density
\begin{align}\label{eq:current}\begin{aligned}
    j^n=\hat{j}\big(\chi_{D_{A^+}}&\sin(\gamma^n+\beta)-\chi_{D_{C^-}}\sin(\gamma^n+\beta-\frac{2\pi}{3})\\&+\chi_{D_{B^+}}\sin(\gamma^n+\beta-\frac{4\pi}{3})\big),
\end{aligned}
\end{align}
with $\hat{j}=I\frac{|D_{A^+}|}{N_w}$, current angle $\beta$ and electrical angle $\gamma^n=N_{\mathrm{pp}}\frac{\pi n}{12 N}.$

This system is discretized by lowest order finite elements on a mesh with 3674 nodes and solved in one monolithic system to incorporate the temporal periodicity. For every position, the rotation of the rotor is incorporated using the harmonic mortar approach \cite{Egger_Mortar}, which also provides a torque formula
\begin{align}\label{eq:torque}
    T^n(a^n)=2N_{\mathrm{pp}}r_{\Gamma}\int_\Gamma\lambda^n(\cur a^n\cdot n_{\Gamma})\circ\rho^n\mathrm{d}s,
\end{align}
where $\lambda^n=\nu_0\cur a^n\cdot n_\Gamma^\perp$ is a Lagrange multiplier and $\rho^n$ realizes the rotational coordinate transformation. Since $\cur u\cdot n_\Gamma=b\cdot n_\Gamma$ is the normal component of the magnetic flux density, this formula shows strong similarities with the well-known torque formula based on the Maxwell stress tensor $r_\Gamma\int_\Gamma h\cdot n_\Gamma^\perp b\cdot n_\Gamma\mathrm{d}s$.
The corresponding Joule losses are computed by 
\begin{align}\label{eq:joule_losses}
    P^J=R_SI^2/2
\end{align} with stator resistance $R_S$ and amplitude value of the current $I$. AC effects on conductors \cite{Hajji2024} are assumed to be included in $R_S$ and are not modeled in detail.
\subsection{Thermal model}
We aim to compute the temperature distribution in the rotor due to EC losses. Since the latency of the heat flux is way higher than the electromagnetic one, we take the average of the losses over one electrical period and solve a single static heat equation. According to \cite{Hameyer_Eddy}, the average EC loss density can be evaluated for every magnet $\Omega_{m_i}, i=1,2,$ by
\begin{align}\label{eq:EC_density}\begin{aligned}
&p^{\mathrm{EC}}_{m_i}(a^0,...,a^{N-1})=\\&\frac{\sigma_m}{N\tau^2}\sum_{n=0}^{N-1}\big(a^n-a^{n-1}-\frac{1}{|\Omega_{m_i}|}\int_{\Omega_{m_i}}a^n-a^{n-1}\mathrm{d}x\big)^2,
\end{aligned}
\end{align}
with $a^0,...,a^{N-1}$ the solution of \eqref{eq:MQS}. Using this, the temperature distribution $\vartheta$ is computed by
\begin{align}\label{eq:heat}\begin{aligned}
    -\mathrm{div}\lambda_\Omega\nabla\vartheta&=p^{\mathrm{EC}}_\Omega(a^0,...,a^{N-1})&&\text{ in } D_R\\
    \lambda_\Omega\nabla\vartheta\cdot n_{\Gamma_{\mathrm{SH}}}&=\beta_{\mathrm{SH}}(\vartheta_0-\vartheta)&&\text{ on }\Gamma_{\mathrm{SH}}\\
    \lambda_\Omega\nabla\vartheta\cdot n_{\Gamma_{\mathrm{AG}}}&=\beta_{\mathrm{AG}}(\vartheta_0-\vartheta)&&\text{ on }\Gamma_{\mathrm{R}},\\
\end{aligned}
\end{align}
where the latter two equations describe the heat flux over the boundaries $\Gamma_{\mathrm{SH}},\Gamma_{\mathrm{AG}}$ to the shaft and the airgap, respectively. Note, that the thermal conductivity $\lambda_\Omega$ is material dependent as well as the source term $p^{\mathrm{EC}}_\Omega=\Chi_{\Omega_{m_1}}p^{\mathrm{EC}}_{m_1}+\Chi_{\Omega_{m_2}}p^{\mathrm{EC}}_{m_2}$. The overall EC losses are computed by
\begin{align}\label{eq:EC_losses}
    P^{\mathrm{EC}}=\ell_z\int_{\Omega_{m_1}\cup\Omega_{m_2}}p^{\mathrm{EC}}_\Omega(a^0,...,a^{N-1})\mathrm{d}x.
\end{align}
\subsection{Structural mechanics}
Since centrifugal forces are increasing with the square of the rotational speed $\omega$, one must not neglect mechanical effects at high speeds. Especially in the fixed iron ring, close to the airgap, stresses might get too big, yielding mechanical failure. As it is common practice in topology optimization, we solve the corresponding linear elastic problem with an artificial weak material in air regions to avoid ill-posed problems. Additionally, we simulate the PMs as "heavy air", i.e., a material with high density but low Young modulus, since the mechanical contact between iron and PMs is not stable to tension, as it was also done in \cite{Lee_MTPAMech}:
\begin{align}
\begin{aligned}\label{eq:elasticity}
    -\mathrm{div}\,\sigma_\Omega(u)&=\rho_\Omega \omega^2x&&\text{ in }D\cup D_{\mathrm{RI}}\\
    \sigma_\Omega(u)n_{\Gamma_R}&=0&&\text{ on }\Gamma_R\\
    u|_{\Gamma_{\mathrm{SH}}}=0,\quad&
    R_{\frac{\pi}{4}}u|_{\Gamma_1}=u|_{\Gamma_2},
\end{aligned}
\end{align}
where $R_\varphi=((\cos\varphi,\sin\varphi)^T,(-\sin\varphi,\cos\varphi)^T)$ is the rotation matrix to realize the periodicity condition on the radial boundaries. We are in the plane stress regime with material dependent Lamé parameters given in Table \ref{tab:data}.
In the optimization, we will impose a constraint on the Von-Mises stress to obtain mechanically feasible designs.

\section{Drive cycle analysis}\label{sec:drivecycle}
\begin{figure}
    \centering
    \includegraphics[width=0.8\linewidth]{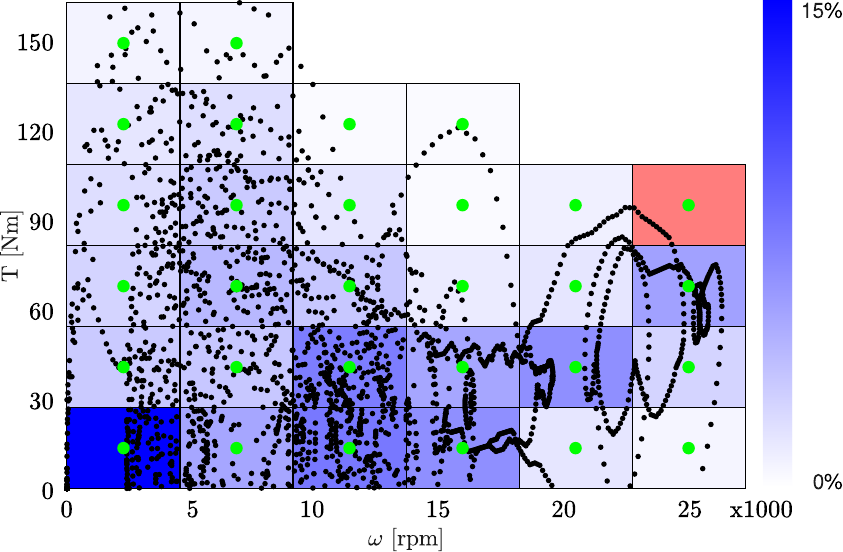}
    \caption{WLTP3 drive cycle. Relative active time $t_k$ as a heat map. Representative OPs $(T_k,\omega_k)$ of selected cells in green. OP considered in thermal analysis in red.}
    \label{fig:drivecycle}
\end{figure}
\begin{figure}
    \centering
    \includegraphics[width=0.48\textwidth]{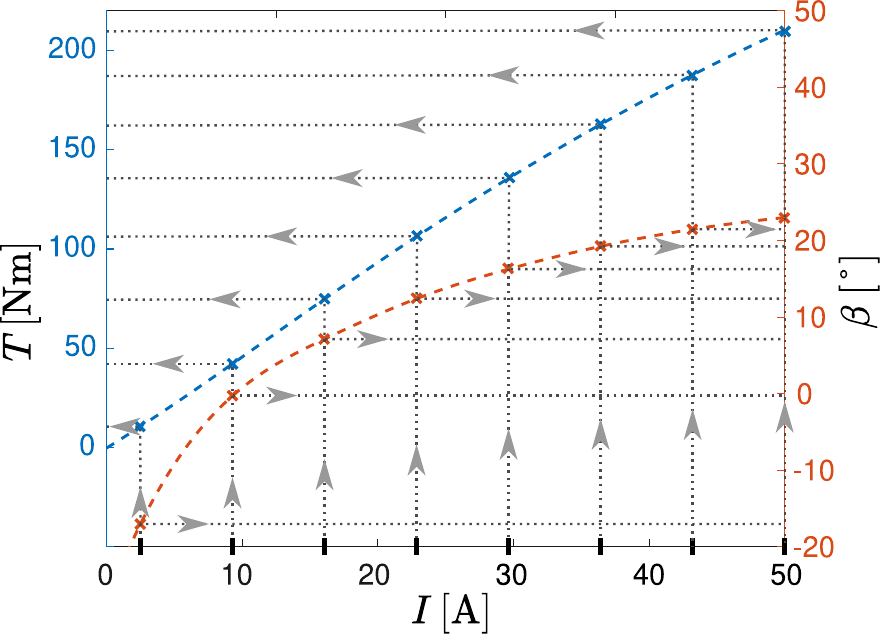}
    \includegraphics[width=0.48\textwidth]{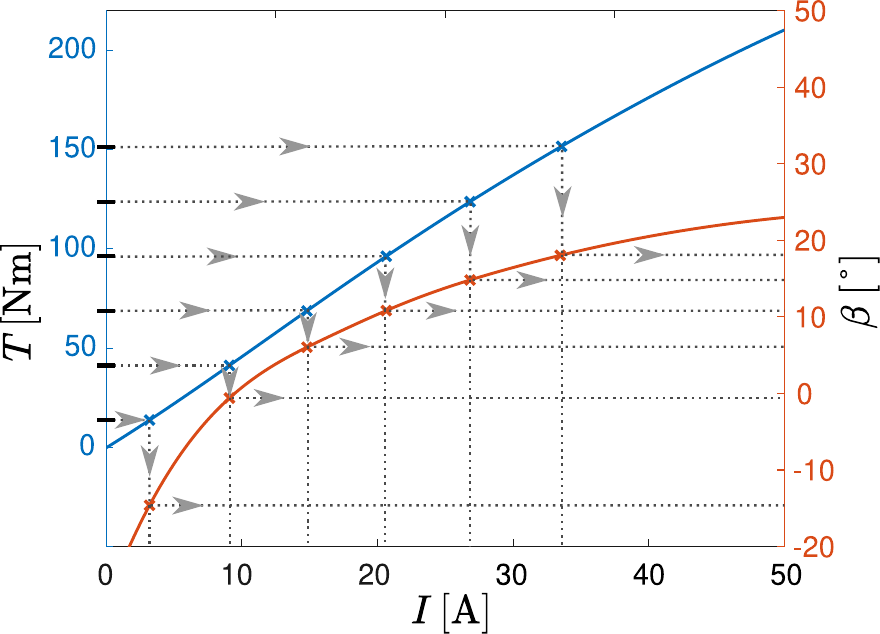}
    \caption{Visualization of the MTPA solution strategy. 
    }
    \label{fig:OP}
\end{figure}

We consider the positive part of the WLTP3 drive cycle \cite{drivecycle}, i.e. the motoric part, scaled to a maximal rotational speed $\omega=27000\mathrm{rpm}$ and a maximal torque $T=160\mathrm{Nm}$, presented in Figure \ref{fig:drivecycle}. We cluster the drive cycle in 30 equally sized cells, represented by their central OP $(\omega_k,T_k)$, and compute the relative active time $t_k$ for each of them. Then, the efficiency \eqref{eq:efficiency} can be approximated by
\begin{align}\label{eq:efficiency_k}
    \E(\Omega)=\frac{\sum_{k=1}^Kt_kP^m_k(\Omega)}{\sum_{k=1}^Kt_k(P_k^m(\Omega)+P_k^J(\Omega)+P_k^{\mathrm{EC}}(\Omega))}.
\end{align}
The mechanical power is the product of speed and torque, which are given by the OPs and therefore independent of the current design $P^m_k(\Omega)=\omega_kT_k$. To compute the losses, we need to find the parameters $I,\beta$ of the excitation current \eqref{eq:current} delivering the desired torque $T_k$. In order to minimize the losses we do this by solving the maximal torque per ampere (MTPA) problem for a given design $\Omega$
\begin{align}\label{eq:mtpa}
        (I_k(\Omega),\beta_k(\Omega))=\underset{(I,\beta)}{\text{argmin}} \,I \text{ s.t. }\overline{\T}(\Omega,(I,\beta))=T_k,
\end{align}
where $\overline{\T}(\Omega,(I,\beta))=\frac{1}{N}\sum_{n=0}^{N-1}T^n(a^n)$ is the average torque \eqref{eq:torque} based on the solution of \eqref{eq:MQS} for a material distribution $\Omega$ and current parameters $(I,\beta)$.
In practice, we take fixed samples of the current amplitude $I_\ell$ in a reasonable range and maximize the average torque
\begin{align}\label{eq:beta}
    \beta_\ell(\Omega)=\underset{\beta\in[0,2\pi)}{\text{argmax}}\,\overline{\T}(\Omega,(I_\ell,\beta))
\end{align}
to find the optimal current angle $\beta_\ell$ based on the system \eqref{eq:MQS} using a gradient descent method, see e.g. \cite{Wiesheu2024}. The corresponding torque value is $T_\ell=\overline{\T}(\Omega,(I_\ell,\beta_\ell))$. Interpolation and inversion of this relation yields a function which maps $T\mapsto I(\Omega)$ according to \eqref{eq:mtpa}. This procedure is sketched in Figure \ref{fig:OP}: On the left we see the samples $I_\ell$ marked on the bottom line. Solving \eqref{eq:beta} yields $\beta_\ell,T_\ell$. On the right we use the interpolation of the points obtained to map $T_k$ first to $I_k$ and then to $\beta_k$.

Since the influence of the ECs on the torque are negligible, we use $\sigma_m=0$ in \eqref{eq:MQS} to speedup the MTPA problem. We obtain, based on \eqref{eq:joule_losses}, the design-dependent Joule losses $P_k^J(\Omega)=R_SI_k(\Omega)^2/2$. Further we compute the EC losses $P^{\mathrm{EC}}_k(\Omega)$ by evaluating \eqref{eq:EC_losses} with $a^0_k(\Omega),...,a^{N-1}_k(\Omega)$ the solution of \eqref{eq:MQS} using $(I_k(\Omega),\beta_k(\Omega))$, now with $\sigma_m>0$.

\section{Topology optimization framework}\label{sec:topop}
\subsection{Topological derivative}
We use the topological derivative $\id^{i\rightarrow j}\J(\Omega)$, a sensitivity information of a shape function $\J$ subject to pointwise changes from material $i$ to material $j$, to update our design. This concept was introduced in \cite{Eschenauer1994-rq} in the context of structural mechanics and first applied to electric machines in \cite{Gangl_Amstutz}. For the electro-thermal coupled problem we rely on \cite{Krenn2025}, which is based on the framework provided in \cite{Gangl_Automated}. In every material point one has to compute the topological derivative for all possible changes, e.g. in a point of a PM $z\in\Omega_{m_1}$ for changes to iron, air and the other magnet. We denote all these sensitivities by a vector-valued quantity
\begin{align}\label{eq:TD}
    \id^i\J(\Omega)(z)=(\id^{i\rightarrow j}\J(\Omega)(z))_{j\in\I, j\neq i}\in \R^{M-1},
\end{align}
for $z\in \Omega_i$, where $\I=\{f,m_1,m_2,a\}$ is the ordered set of $M=|\I|=4$ materials present in the design space $D$. The evaluation of the topological derivative is based on the solutions of the involved partial differential equations and the corresponding adjoint equations. For the sake of brevity, we do not state them; for all problems considered here, they can be found in \cite{Krenn2025}.
\subsection{Level set framework}
Following the work of \cite{Amstutz_Levelset} and the extension to multiple materials \cite{Gangl_Multi}, we use a vector-valued level set function $\psi:D\rightarrow\R^{M-1}$ to represent the design by
\begin{align}
    z\in \Omega_i\Leftrightarrow\psi(z)\in S_i, i\in\I,
\end{align}
with sectors $S_i=\{x\in\R^{M-1}:|V_i-x|<|V_j-x|, j\in\I,j\neq i\}$, where $V_i$ are the $M$ vertices of the $M-1$ dimensional regular unit simplex.
This means the material in each point is defined by the sector, into which the level set function is pointing. To map between topological derivatives and level set functions, we need the matrices 
$N_i=\left((V_i-V_j)^T\right)_{j\in\I,j\neq i}$.
As shown in \cite{Gangl_Multi}, a design is optimal if there exists a positive number $c>0$ such that $N_i\psi(z)=c\id^i\J(\Omega_\psi)(z)$ for all material points $z\in\Omega_i, i\in\I$, where $\Omega_{\psi}$ is the design represented by $\psi$. This motivates the update of the level set function
\begin{align}\label{eq:ls_update}
    \psi^{m+1}=(1-s)\psi^m+s\sum_{i\in\I}\Chi_{\Omega_i}N_i^{-1}\id^i\J(\Omega_{\psi^m}),
\end{align}
where the stepsize $s$ is chosen small enough such that $\J(\Omega_{\psi^{m+1}})<\J(\Omega_{\psi^m})$; for more details see \cite{Gangl_Multi, Krenn2025}.
\subsection{Efficiency maximization}\label{sec:efficiency}
\begin{figure}
    \centering
    \includegraphics[width=0.9\textwidth]{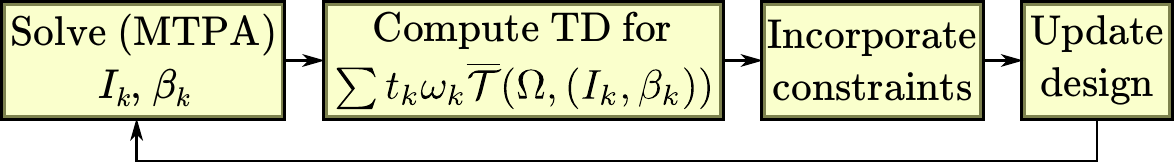}
    \caption{Flowchart of the drive cycle design optimization}
    \label{fig:flow}
\end{figure}

To maximize the efficiency, we use the following workflow, displayed in Figure \ref{fig:flow}: First, we solve the MTPA problem \eqref{eq:mtpa}, to determine the current parameters $(I_k,\beta_k)$ for the actual design $\Omega$. Next we compute the topological derivative \eqref{eq:TD} of the negative weighted average torque 
\begin{align}\label{eq:cost}
    \J(\Omega)=-\sum_{k=1}^Kt_k\omega_k\overline{\T}(\Omega,(I_k,\beta_k))
\end{align}
for these current parameters. We have the minus, since the optimization framework is minimizing $\J(\Omega)$ and we want to maximize the performance. 
After updating the design according to \eqref{eq:ls_update}, we go back to step one. The topological derivative of the average torque functional, we use here, is given in \cite[Eq. 71]{Krenn2025}.
\section{Constraint handling}\label{sec:constraints}
\subsection{Temperature constraint}
Since the EC losses grow quadratically with the rotation speed, the highest temperatures will occur for the OP with first, maximal speed and second, maximal torque. This OP is highlighted in red in Figure \ref{fig:drivecycle}. We aim to control the maximal PM temperature in both electrical and thermal steady states for this OP. To do so, we solve the heat equation \eqref{eq:heat} with the average EC loss density \eqref{eq:EC_density} based on the solutions $a^0_k,...,a_k^{N-1}$ of \eqref{eq:MQS}. A design is now feasible if $\vartheta\le\vartheta^*$ for a given temperature bound $\vartheta^*$. We model this by the functional
\begin{align}
    \C_t(\Omega)=\int_{\Omega_{m_1}\cup\Omega_{m_2}}\left(\max\left\{1,\frac{\vartheta}{\vartheta^*}\right\}-1\right)^2\id x
\end{align}
which is zero, if and only if the constraint is fulfilled. Since it is always positive, by construction, we can add it to our objective \eqref{eq:cost} with some positive weight 
\begin{align}\label{eq:cost_temp}
    \J_t(\Omega)=\J(\Omega)+w_t\C_t(\Omega).
\end{align}
If $w_t$ is chosen big enough, $\C_t(\Omega)$ will be close to zero, in order to minimize $\J_t(\Omega).$ The topological derivative for this constraint functional, based on the electro-thermal problem, can be found in \cite[Eq. 86]{Krenn2025}.
\subsection{Von-Mises stress constraint}
Our second criterion of feasibility is on the squared Von-Mises stress $s_{\mathrm{VM}}=\frac{1}{2}(3|\sigma_\Omega(u)|^2-\mathrm{tr}(\sigma_\Omega(u))^2)\le (\sigma^*)^2$ in the design domain $D$ and the iron ring $D_{\mathrm{RI}}$, with $u$ being the solution of \eqref{eq:elasticity} for the maximal rotation speed $\omega=25000 \mathrm{rpm}$. Similarly as for the temperature constraint, we mimic this constraint by a positive functional
\begin{align*}
    \C_{\mathrm{VM}}(\Omega)=\int_{D\cup D_{\mathrm{RI}}}\left(1+\left(\frac{s_{\mathrm{VM}}}{(\sigma^*)^2}\right)^p\right)^{\frac{1}{p}}\id x.
\end{align*}
We have to choose here a different function in order to assure the topologically differentiability, see \cite{Amstutz_vonMises}. For $p=\infty$ we have again that $\C_{\mathrm{VM}}(\Omega)$ is zero if and only if the constraint is fulfilled. In practice we choose $p=16$ which turned out to be a sufficiently good approximation preserving numerical stability. We add this to the objective with the weight $w_{\mathrm{VM}}$
\begin{align}\label{eq:cost_all}
    \J_{t,\mathrm{VM}}(\Omega)=\J(\Omega)+w_t\C_t(\Omega)+w_{\mathrm{VM}}\C_{\mathrm{VM}}(\Omega).
\end{align}
The topological derivative of $\C_{\mathrm{VM}}(\Omega)$ is taken from \cite[Eq. 94]{Krenn2025}, which is based on \cite{Amstutz_vonMises}.
\section{Numerical optimization results}\label{sec:results}
\begin{figure}
    \centering
    \begin{tabular}{ccc}
    \includegraphics[width=0.48\textwidth,trim=8 8 8 8,clip]{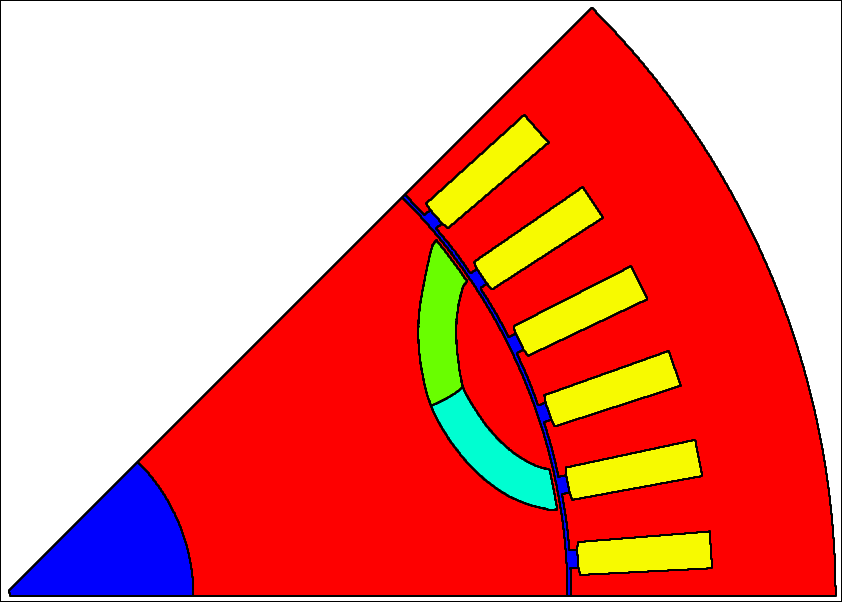}&
    \includegraphics[width=0.16\textwidth,trim=380 8 260 160,clip]{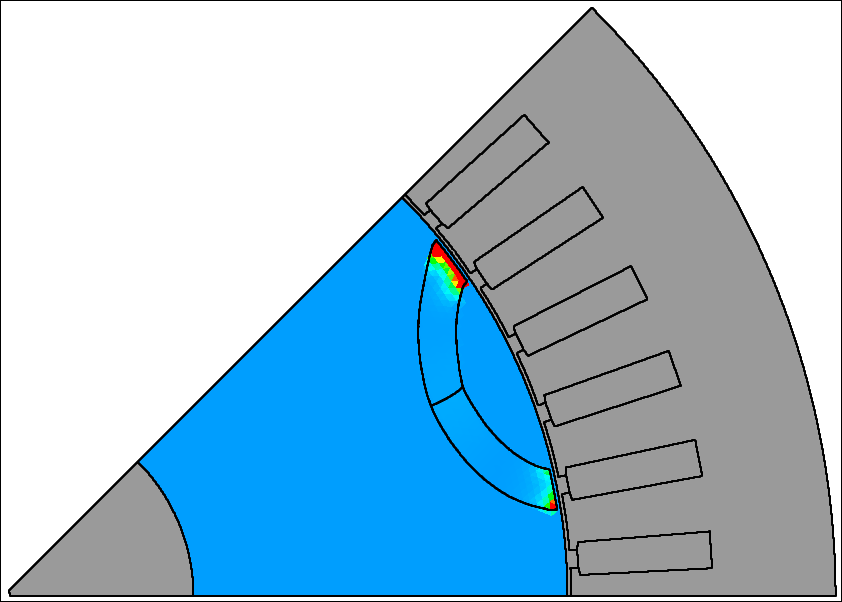}&
    \includegraphics[width=0.16\textwidth,trim=380 8 260 160,clip]{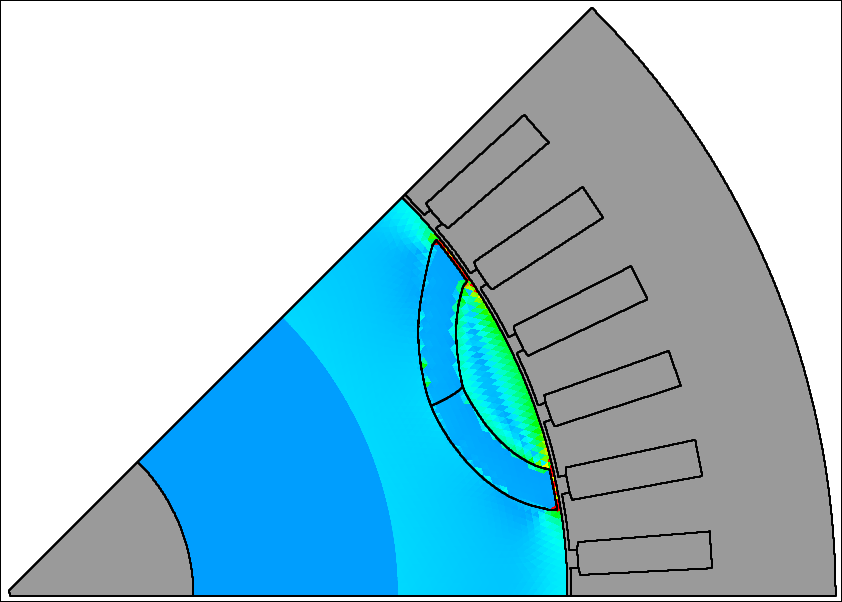}\\
    \includegraphics[width=0.48\textwidth,trim=8 8 8 8,clip]{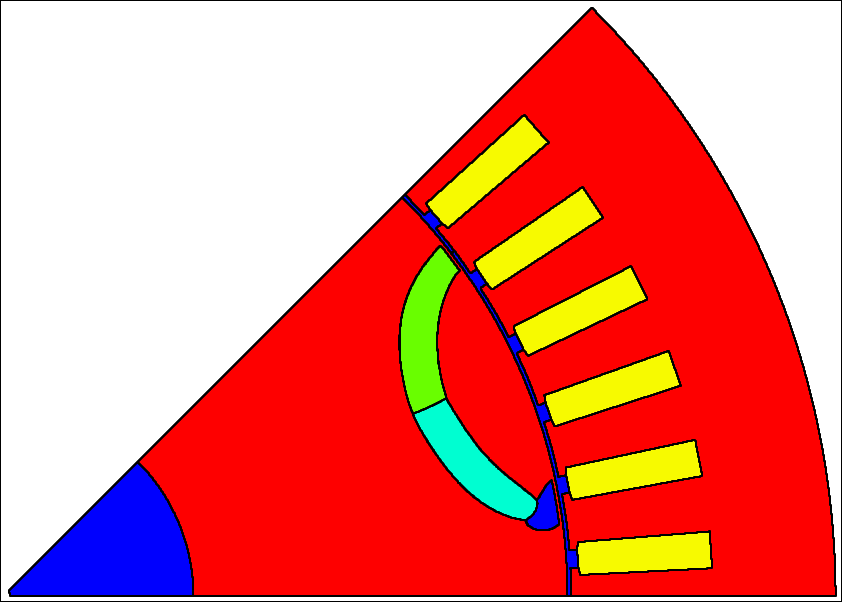}&
    \includegraphics[width=0.16\textwidth,trim=380 8 260 160,clip]{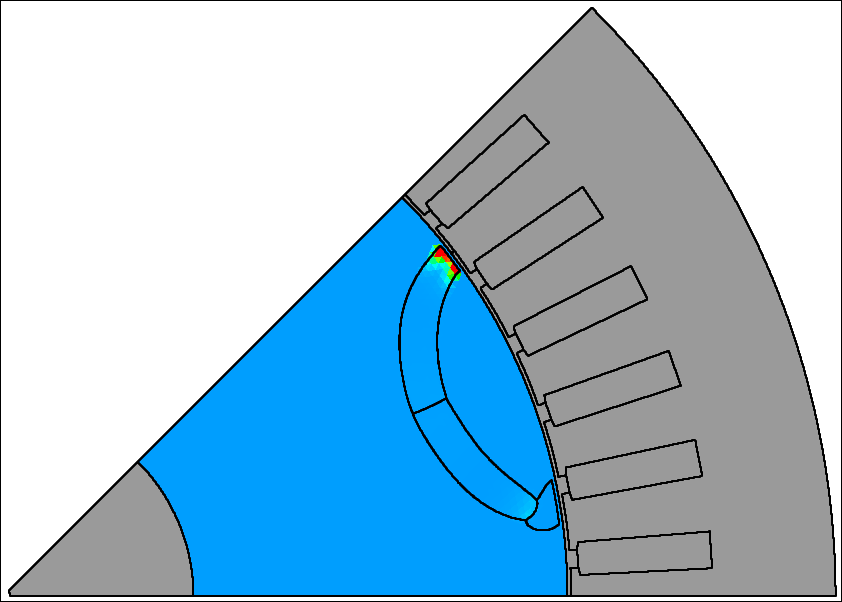}&
    \includegraphics[width=0.16\textwidth,trim=380 8 260 160,clip]{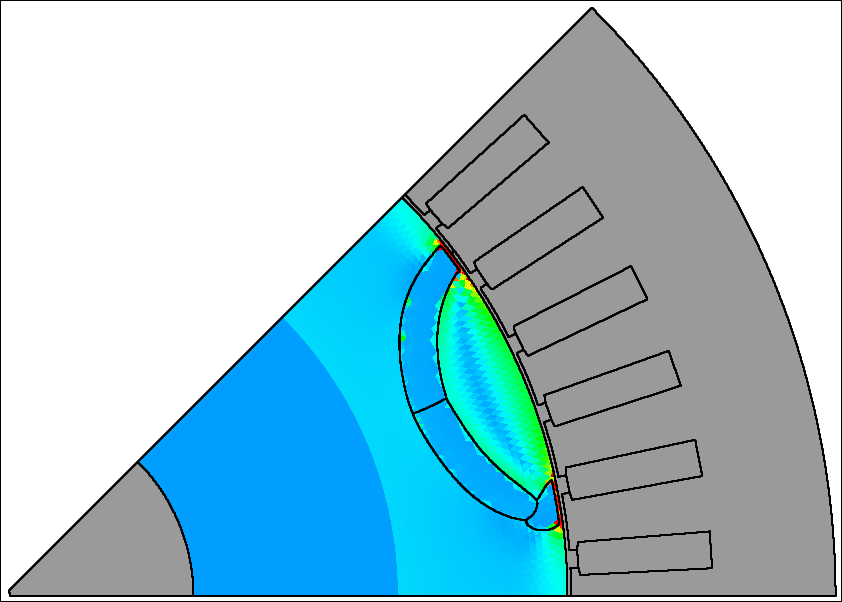}\\
    \includegraphics[width=0.48\textwidth,trim=8 8 8 8,clip]{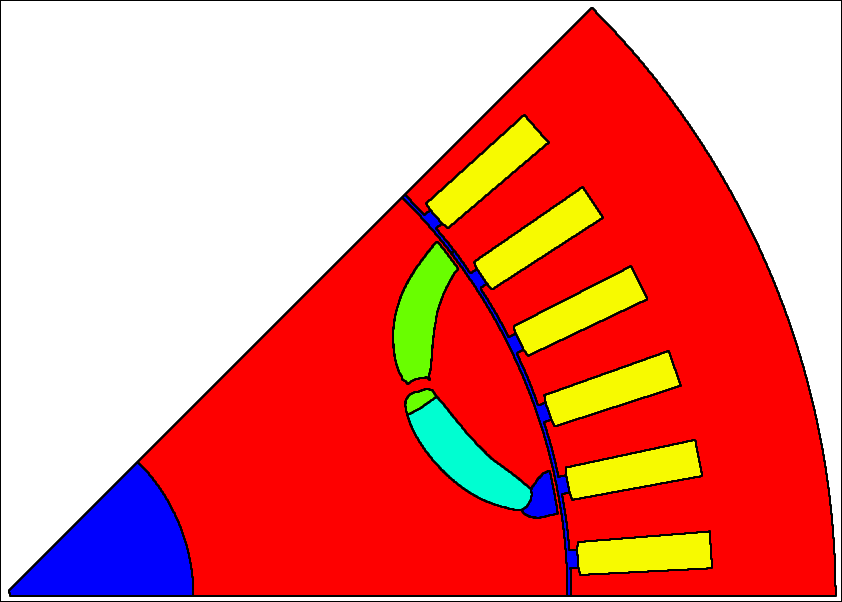}&
    \includegraphics[width=0.16\textwidth,trim=380 8 260 160,clip]{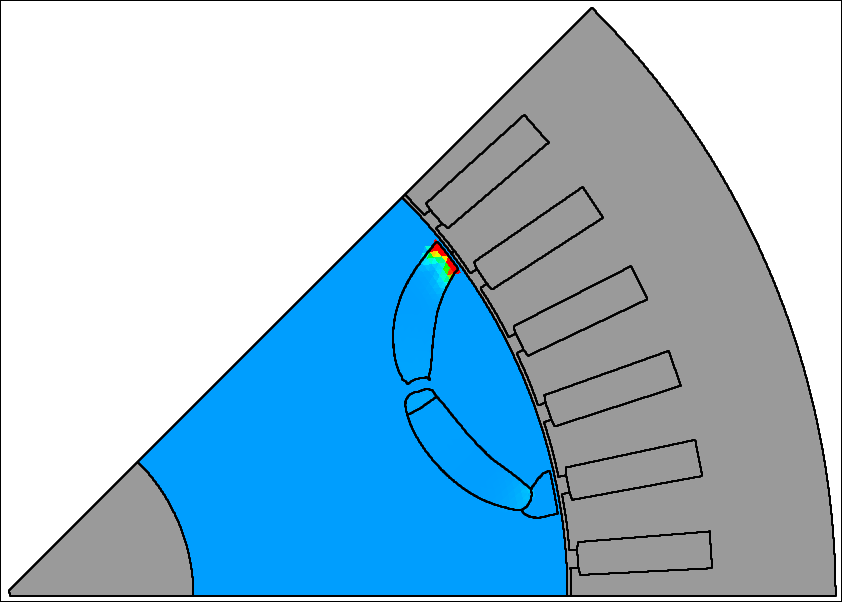}&
    \includegraphics[width=0.16\textwidth,trim=380 8 260 160,clip]{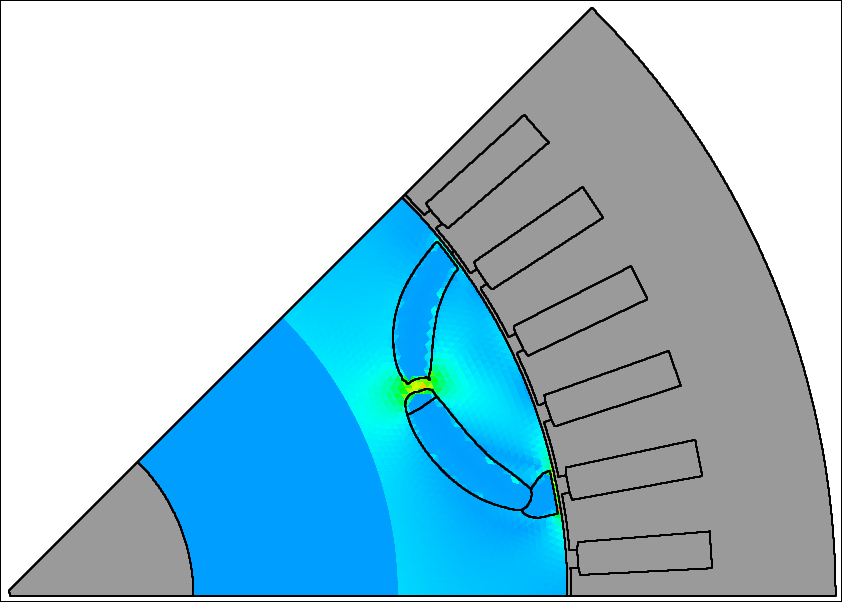}
    \end{tabular}
\caption{Designs (left), EC loss density (mid) and Von-Mises stresses (right) of initial (1. line), unconstrained (2. line), thermally constrained (3. line) and thermally and mechanically constrained optimization (4. line).}    \label{fig:all}
\end{figure}

First, we start by analyzing a state-of the art design $\Omega_{\mathrm{ini}}$, shown in  Figure \ref{fig:machine}. Temperature and stresses are in the feasible range, as noted in Table \ref{tab:results}. Starting from this design, we run optimizations considering different constraints. The obtained designs are shown in Figure \ref{fig:all} together with the EC loss density (which we preferred to plot, since the temperature variation is very small for every design) and the Von-Mises stress distribution.

By the procedure as described in Section \eqref{sec:efficiency}, we can increase the efficiency defined in \eqref{eq:efficiency_k} from $\E(\Omega_{\mathrm{ini}})=95.74\%$ to $\E(\Omega^*)=97.34\%$. This comes with the drawback of violating both, the thermal constraint $\vartheta^*=80^\circ\mathrm{C}$ by $31^\circ$ and the Von-Mises stress constraint $\sigma^*=500\mathrm{MPa}$ by $800\mathrm{MPa}$.

We are able to reduce the temperature below the desired $\vartheta^*=80^\circ\mathrm{C}$ by considering $\J_t(\Omega)$ \eqref{eq:cost_temp} with $w_t=10^7$ in the optimization while almost maintaining the improved efficiency $\E(\Omega_t^*)=97.26\%$. In the plot of the EC loss densities, we can see, comparing second and third line, that the areas of high losses, displayed in red, get smaller. Still, mechanical stresses are too high, especially in the iron ring $D_{\mathrm{RI}}$.

Finally, we consider also the stress constraints in the optimization by $\J_{t,\mathrm{VM}}(\Omega)$ \eqref{eq:cost_all} with $w_{\mathrm{VM}}=10^{10}$. The resulting design $\Omega_{t,\mathrm{VM}}$ has again a slightly lower efficiency $\E(\Omega_{t,\mathrm{VM}}^*)=97.00\%$, stays inside the temperature bound and fulfills the stress constraint $\sigma^*=500\mathrm{MPa}$. Mechanical stability, as can be seen by the absence of red zones in the stress plot, is achieved by the iron bridge between the PMs.
\begin{table}[]
    \centering
    \begin{tabular}{ccccc}
     Design& $\mathcal{E}(\Omega)$&$\max\vartheta$&$\max \sqrt{s_{\mathrm{VM}}}$&iterations\\
     \hline 
    $\Omega_{\mathrm{ini}}$ &95.47\% & \phantom{1}$51^\circ\mathrm{C}$&$\phantom{1}320\mathrm{MPa}$&-\\
    $\Omega^*$ &97.34\% &$111^\circ\mathrm{C}$&$1300\mathrm{MPa}$&105\\
    $\Omega^*_t$ &97.26\% &\phantom{1}$74^\circ\mathrm{C}$&$1654\mathrm{MPa}$&229\\    
    $\Omega^*_{t,\mathrm{VM}}$ &97.00\% &\phantom{1}$80^\circ\mathrm{C}$ &$\phantom{1}437\mathrm{MPa}$&394  
    \end{tabular}
    \caption{Evaluation of efficiency, maximal temperature and Von-Mises stress together with the iteration number. }
    \label{tab:results}
\end{table}

\section{Conclusion and outlook}\label{sec:conclusion}
In this work we successfully introduced an innovative methodology to handle multiple OPs in a multi-physical, multi-material topology optimization subject to constraints on both, PM temperature and mechanical stresses. The combination of these different aspects is, to our best knowledge, new in the free form design optimization of electric machines. 

There are various options for further work. One can improve the loss models including iron losses or AC winding losses. It would be also interesting to incorporate cooling of the machine throughout the drive cycle considering the temporal heating process.

\section*{Acknowledgment}
The work of P.G., N.K. and S.Sch. is partially supported by the joint DFG/FWF Collaborative Research Centre CREATOR (DFG: Project-ID 492661287/TRR 361; FWF: 10.55776/F90) at TU Darmstadt, TU Graz, JKU Linz and RICAM Linz. P.G. is partially supported by the State of Upper Austria.


\end{document}